
\magnification1200
\input amssym.def
\input amssym.tex 
\def\SetAuthorHead#1{}
\def\SetTitleHead#1{}
\def\NoindentAfter{\everypar={\setbox0=\lastbox\everypar={}}}
\def\H#1\par#2\par{{\baselineskip=15pt\parindent=0pt\parskip=0pt
 \leftskip= 0pt plus.2\hsize\rightskip=0pt plus.2\hsize
 \bf#1\unskip\break\vskip 4pt\rm#2\unskip\break\hrule
 \vskip40pt plus4pt minus4pt}\NoindentAfter}
\def\HH#1\par{{\bigbreak\noindent\bf#1\medbreak}\NoindentAfter}
\def\HHH#1\par{{\bigbreak\noindent\bf#1\unskip.\kern.4em}}
\def\th#1\par{\medbreak\noindent{\bf#1\unskip.\kern.4em}\it}
\def\endth{\medbreak\rm}
\def\pf#1\par{\medbreak\noindent{\it#1\unskip.\kern.4em}}
\def\df#1\par{\medbreak\noindent{\it#1\unskip.\kern.4em}}

\let\Roster\bgroup\let\endRoster\egroup
\def\\{}\def\text#1{\hbox{\rm #1}}
\def\mop#1{\mathop{\rm\vphantom{x}#1}\nolimits}
\def\MaxReferenceTag#1{}
\def\qedbox{\vrule width2mm height2mm\hglue1mm\relax}
\def\qed{\ifmmode\qedbox\else\hglue5mm\unskip\hfill\qedbox\medbreak\fi\rm}

\def\cite#1{{\bf[#1]}}
\def\Em#1{{\it #1\/}}
\def\Bib#1\par{\bigbreak\bgroup\centerline{#1}\medbreak\parindent30pt
 \parskip2pt\frenchspacing\par}
\def\endBib{\par\egroup}
\newdimen\Overhang
\def\rf#1{\par\noindent\hangafter1\hangindent=.7 true in
     \setbox0=\hbox{[#1]}\Overhang\wd0\advance\Overhang.4em\relax
     \ifdim\Overhang>\hangindent\else\Overhang\hangindent\fi
     \hbox to \Overhang{\box0\hss}\ignorespaces}

\def\bbZ{{\Bbb Z}}
\def\Coordinates{\bigbreak\bgroup\parindent=0pt\obeylines}
\def\endCoordinates{\egroup}

\def\isom{\mathop{\rm Isom}}
\def\sol{{\rm\bf Sol}}

\tolerance 10000
\def\denom{\mop {denom}}
\SetTitleHead{}
\SetAuthorHead{}

\H Solvable Baumslag-Solitar Groups Are Not Almost Convex

Charles F.  Miller III and Michael Shapiro\footnote{*}{We wish to
thank the ARC for support.}

{\narrower{The arguments of Cannon, Floyd, Grayson and Thurston
\cite{CFGT} showing that solvegeometry groups are not almost convex
apply to solvable Baumslag-Solitar groups.\par}}

\HH Introduction

The property of {almost convexity} was first introduced by Cannon in
\cite{C}.  This property has very geometric in flavor, being defined in
terms of the geometry of the Cayley graph.  If the Cayley graph is
almost convex then there are efficient algorithms for calculating in
$G$, or, if you like, constructing the Cayley graph of $G$.

Given a group $G$ and a finite generating set $C \subset G$, the
\Em{Cayley graph} of $G$ with respect to $C$ is the directed,
labeled graph whose vertices are the elements of $G$ and whose
directed edges are the triples $(g,c,g')$ such that $g,g'\in G$, $c
\in C$ and $g'=gc$.  Such an edge is directed from $g$ to $g'$ and is
labeled by $c$.  We denote this Cayley graph by
$\Gamma=\Gamma_{C}(G)$.  We will assume that $C$ is closed under
inverses.

A Cayley graph has a natural base point $1 \in G$ and a natural path
metric $d = d_{C}$ which results from identifying each edge with the
unit interval.  Each element of $G$ has a natural \Em{length}
$\ell(G)=\ell_{C}(g)=d_{C}(1,g)$.  We define the \Em{ball of radius}
$n$ to be
$$B(n) = \left\{ x \in \Gamma \mid d_{C}(1,x) \le n \right\}.$$
The group $G$ is \Em{almost convex $(k)$} with respect to $C$ if there
is $N=N(k)$ so that if $g,g' \in B(n)$ and $d_{C}(g,g')\le k$ then
there is a path $p$ from $g$ to $g'$ inside $B(n)$ whose length is at
most $N$.  $G$ is \Em{almost convex} with respect to $C$ if it is almost
convex $(k)$ with respect to $C$ for each $k$.  It is a theorem of Cannon
\cite{C} that if $G$ is almost convex $(2)$ with respect to $C$ then
$G$ is almost convex with respect to $C$.

Thiel \cite{T} has show that almost convexity is not a group property,
i.e, that there are groups which are almost convex with respect to one
generating set but not another.  Almost convexity is fairly well
understood for the fundamental groups of closed $3$-manifold groups
with uniform geometries \cite{SS}.  The solvegeometry case was covered
in a beautiful paper of Cannon, Floyd, Grayson, and Thurston
\cite{CFGT}{\parindent 0 true in\footnote{$^{\rm\dag}$} {Cannon has pointed
out to us that there are problems with some of the details in their
paper.  These concern the relationship between lengths in the given
group $G\subset \isom(\sol)$ and the finite index subgroup
$A=G\cap\sol$.  These problems are easily fixed and in our view, their
paper remains quite beautiful.}}.  They show that any group which acts
co-compactly, discretely by isometries on \sol\  cannot be almost convex
with respect to any generating set.  In this paper we show that their
arguments apply to solvable Baumslag-Solitar groups
$$G=B_{1,p}=\langle a,t \mid t^{-1}at=a^{p}\rangle$$ 
with $|p| > 1$.

\th Theorem

Let $G= B_{1,p}$ be a solvable Baumslag-Solitar group with $|p| > 
1$.  Then $G$ is not almost convex with respect to any generating 
set.\endth

\HH Proof of the Theorem

Let $G=B_{1,p}$ with $|p| > 1$.  Then $G$ has the form
$$G=\bbZ\left[{1/ p}\right]\rtimes\bbZ,$$
where the generator of $\bbZ$ acts via multiplication by $p$.  Thus
each element of $G$ has the form $(f,c)$ where $f$ is a fraction of
the form $$f={m \over p^{n}},$$ and $m,n,c\in \bbZ$.  For each
element $(f,0) \in \bbZ[1/p] \subset G$, we will take $|(f,0)|=|f|$
and if $n$ is minimal such that $f=m/p^{n}$, we will say that $|p|^{n}$
is the \Em{denominator} of $(f,0)$ written $\denom(f,0)$.

We fix a generating set
$$C=\{(f_{i},c_{i})\}$$
which we assume is closed under inverses.  We take
$$\eqalign{c&= \max\{c_{i} \mid (f_{i},c_{i}) \in C \}\cr
  f^{*} &= \max\{f_{i} \mid (f_{i},c_{i}) \in C \}\cr
  f^{**} &= \max\{\denom(f_{i},0) \mid (f_{i},c_{i}) \in C \}.}$$
We assume $(f_{*},c)\in C$ realizes the first of these maxima.
Notice that $c>0$.

We need the following two lemmas which give information about distance with
respect to the generating set $C$.

\th Lemma 1

There is a constant $M$ so that if $(f,0) \in B(n)$ then either $|f|
\le  M |p|^{nc\over 4}$ or $\denom(f) \le M |p|^{nc\over 4}$.  Further, both $|f|
\le  M |p|^{nc\over 2}$ and $\denom(f) \le M |p|^{nc\over 2}$
\endth

\pf Proof

First observe the following product formula:
$$(f_{1},c_{1}) \ldots (f_{n},c_{n}) = \left( \sum_{i=1}^{n} f_{i}
p^{(0-c_{1}-\cdots - c_{i-1})}, \sum_{i=1}^{n} c_{i}\right).$$
Since $(f,0) \in B(n)$, $(f,0)$ can be written as such a product where
each $(f_{i},c_{i})$ is in $C$ and $ \sum_{i=1}^{n} c_{i}=0$.
For each $i$, $i=1,\ldots,n$, we set $e_{i}=0-c_{1}-\dots - c_{i-1}$.
Then $e_{i}$ is positive for at most $n/ 2$ values of $i$, or
$e_{i}$ is negative for at most $n/ 2$ values of $i$.  Suppose
that $e_{i}$ is positive for at most $n/ 2$ values of $i$.  We
then have
$$\eqalign{|f| &= \sum_{i=1}^{n} f_{i}p^{(0-c_{1}-\cdots - c_{i-1})}\cr
&\le  f^{*} \sum_{i=1}^{n} |p|^{(0-c_{1}-\cdots -
c_{i-1})}\cr
&=  f^{*}\left(\sum_{e_{i}\le 0} |p|^{e_{i}} + \sum_{e_{i}> 0}
|p|^{e_{i}}\right) \cr
&\le f^{*} \left(n + \sum_{e_{i}> 0} |p|^{e_{i}}\right).}$$
Let us enumerate the $\{e_{i} \mid e_{i} >0 \}$ as $i$ increases, so
that these are the $\lfloor{n/ 2}\rfloor$-tuple $(e_{i_{1}},
e_{i_{2}}, \dots, e_{i_{\lfloor {n\over 2}\rfloor}})$.  If there are
less than $\lfloor{n/2}\rfloor$ of these, we will consider any final
entries in this list to be $0$.  We now take the $m$-tuple
$(e'_{i_{1}}, e'_{i_{2}}, \dots, e'_{i_{m}})$ to be
$(c,2c,\dots,{nc\over 4},{nc\over 4}, \dots, 2c,c)$.  Here $m$ is
either $\lfloor n/2 \rfloor$ or $\lfloor n/2 \rfloor+1$. It is not hard
to see that for each $j$, $e_{i_{j}}\le e'_{i_{j}}$.  Consequently,

$$|f| \le f^{*} \left( n + 2|p|^{c}+ 2|p|^{2c} +\cdots+2|p|^{nc\over 4}
\right),$$
and thus for suitable choice of $M'$, $|f|\le M'|p|^{nc\over4}$.

On the other hand if more than $n/ 2$ of the $e_{i}$ are positive,
then less than $n/ 2$ of them are negative, and in particular, the
most negative any of these can be is ${-nc\over 4}$.  It immediately
follows that
$$\denom(f,0) \le f^{**}|p|^{{nc\over 4}}.$$

Taking $M=\max\{M',f^{**}\}$ completes the proof of the first part of Lemma
1.  After suitably enlarging $M$, a completely similar proof gives the
simultaneous bound on $|f|$ and $\denom(f)$\qed

We also need the following observation.

\th Lemma 2

If $h,h' \in \bbZ[1/p] \subset G$ with $d_{C}(h,h') \le r $ then
$||h|-|h'||\le M |p|^{rc\over 2} $ and $|\denom(h) - \denom(h')| \le M
|p|^{rc\over 2} $.\endth

\pf Proof

If $d_{C}(h,h') \le r$ then (using additive notation in $\bbZ[1/p]$)
$h - h' \in B(r)$.  Thus
$$||h|-|h'|| \le |h-h'| \le M |p|^{rc\over 2} .$$

On the other hand $h = h' +(h-h')$ and $h'= h -(h-h')$ so we have
$$\denom(h) \le\max\{\denom(h'), \denom(h-h')\} \le \denom(h') +
\denom(h-h'),$$
and
$$\denom(h') \le\max\{\denom(h), \denom(h-h')\} \le \denom(h) +
\denom(h-h').$$
Consequently,
$$|\denom(h) - \denom(h')| \le \denom(h-h') \le M |p|^{rc\over 2} .$$\qed

We now return to the proof of the Theorem.  For each $k>0$ we take
$$\eqalign{T_k &=(f_{*},c)^{-k} (1,0) (f_{*},c)^{k} = (p^{ck},0) \cr
           S_k &=(f_{*},c)^{k} (1,0) (f_{*},c)^{-k} = (p^{-ck},0) }$$
We then have $T_{k}S_{k}=S_{k}T_{k}$.  For some $j$ which we will fix
later, we take
$$\eqalign{\alpha_{k} &= S_{k}T_{k}(f_{*},c)^{-j} \cr
           \beta_{k}  &= T_{k}S_{k}(f_{*},c)^{j}}$$
If we take $\ell=\ell_{C}(1,0)$ and $k>j$, then $\alpha_{k}$ and
$\beta_{k}$ both lie in $B(4k+2\ell-j)$ and within distance $2j$ of
each other.

Suppose that, contrary to hypothesis,  $G$ is almost convex.
Then there is a constant $N=N(2j)$ so that $\alpha_{k}$ and $\beta_{k}$
are joined by a path of length at most $N$ lying entirely within
$B(4k+2\ell-j)$.  The second coordinates of points along this path
vary from $-jc$ to $+jc$ changing by at most $\pm c$ along each edge.
In particular, this path must pass through a point
$P'_{k}$ of the form $(g_{k},i)$ with $|i|\le {c/ 2}$.  We take
$$\epsilon = \max\left\{\ell_{C}(0,i) \mid |i| \le {c/ 2}\right\}.$$
It follows that the point
$$P_{k}=(g_{k},0)$$
lies within $B(4k+2\ell+\epsilon-j)$ and within distance ${N/2}+\epsilon+j$ of
$S_{k}T_{k}$.  Notice that the distance from $P_{k}$ to $S_{k}T_{k}$
is bounded by a constant independent of $k$.
It is this fact that we will contradict, thus showing $G$ is not almost
convex in the given generating set.

  From Lemma 1 above it follows that either
$$|P_{k}| \le M |p|^{c({k+{2\ell +\epsilon\over 4}-{j\over 4}})}$$
or
$$\denom(P_{k}) \le M |p|^{c({k+{2\ell +\epsilon\over 4}-{j\over 4}})}.$$
Let us fix $j$ so that
$$|p|^{j\over 4} > M|p|^{{2\ell +\epsilon\over 4c}+1}$$
and hence
$$M |p|^{c({k+{2\ell +\epsilon\over 4}-{j\over 4}})}=
|p|^{kc}(M|p|^{{2\ell +\epsilon\over 4c}}|p|^{-{j\over 4}})
 \le |p|^{kc-1}.$$
It then follows that either
$$|P_{k}| \le |p|^{kc-1}$$
or
$$\denom(P_{k}) \le |p|^{kc-1}.$$

Now $|S_kT_{k}| = |p^{kc} + p^{-kc}|$ and $ \denom(S_{k}T_{k})=|p|^{kc}$.
Hence either
$$|S_kT_{k}|- |P_{k}|
\ge |p^{kc} + p^{-kc}|-|p|^{kc-1} = |p|^{kc-1}(|p +  p^{-2kc+1}|-1)>|p|^{kc-1}$$
or
$$\denom(S_kT_{k}) - \denom(P_k)
\ge |p|^{kc}-|p|^{kc-1} = |p|^{kc-1}(|p|-1)\ge |p|^{kc-1}.$$
In either case by Lemma 2,
as $k\to\infty$ the distance $d_C(P_k,S_kT_k)$
increases without bound.
But our assumption of almost convexity
implied $d_C(P_k,S_kT_k) \le N/2 +\epsilon + j$
which is a constant. This is a contradiction.
Hence $G$ is not almost convex.\qed

\Bib{References}
\MaxReferenceTag{CFGT}

\rf{C} J. Cannon, Almost convex groups, Geom. Dedicata 22 197--210 (1987).

\rf{CFGT} J. Cannon, W. Floyd, M. Grayson and W. Thurston, Solvgroups
are not almost convex, Geom. Dedicata 31 no.\ 3 292--300 (1989).

\rf{SS}  M.~Shapiro and M.~Stein,  Almost convexity and
the eight geometries,  Geom. Dedicata, 55, 125--140, (1995).

\rf {T} C. Theil, Zur fast-konvexit\"at einiger nilpotenter gruppen,
Doctoral dissertation, Bonn, 1991.

\endBib

\bigskip

\Coordinates

Department of Mathematics and Statistics
University of Melbourne
Parkville, VIC 3052
Australia

\endCoordinates

\bye